\documentclass{amsart}

\usepackage{amssymb,amsmath,latexsym,amsthm}

\newcommand{\norm}[1]{\ensuremath{\left\Vert #1 \right\Vert}}
\newcommand{\abs}[1]{\ensuremath{\left\vert #1 \right\vert}}
\newcommand{\intdist}[1]{\ensuremath{\left\vert\left\langle #1
      \right\rangle\right\vert}}
\newcommand{\infnorm}[1]{\ensuremath{\left\Vert #1 \right\Vert_\infty}}
\newcommand{\hdim}{\text{dim}_{\text{H}}}
\newcommand{\matA}[1]{\widetilde{#1}}
\newcommand{\column}[2]{\ensuremath{#2^{(#1)}}}
\DeclareMathOperator{\windim}{windim}

\newtheorem{theorem}{Theorem}[section]
\newtheorem{lemma}[theorem]{Lemma}

\newtheorem{corollary}[theorem]{Corollary}

\theoremstyle{definition}
\newtheorem*{definition}{Definition}

\begin{document}

\title[Badly approximable systems of linear forms]{Badly approximable
  systems of linear forms over a field of formal series}

\author{Simon Kristensen}
\address{Simon Kristensen\\
  School of Mathematics \\
  The University of Edinburgh \\ 
  James Clerk Maxwell Building \\
  Kings Buildings \\
  Mayfield Road \\
  Edinburgh \\
  EH9 3JZ \\
  Scotland}
\email{Simon.Kristensen@ed.ac.uk}

\begin{abstract}
  We prove that the Hausdorff dimension of the set of badly
  approximable systems of $m$ linear forms in $n$ variables over the
  field of Laurent series with coefficients from a finite field is
  maximal.  This is a analogue of Schmidt's multi-dimensional
  generalisation of Jarn\'ik's Theorem on badly approximable numbers.
\end{abstract}

\maketitle

\section{Introduction}

Let $\mathbb{F}$ denote the finite field of $k = p^r$ elements, where
$p$ is a prime and $r$ is a positive integer. We define
\begin{equation}
  \label{eq:1}
  \mathcal{L} = \left\{\sum_{i=-n}^\infty a_{-i} X^{-i} : n \in
    \mathbb{Z}, a_i \in \mathbb{F}, a_n \neq 0\right\} \cup \{0\}.
\end{equation}
Under usual addition and multiplication, this set is a field,
sometimes called \emph{the field of formal Laurent series with
  coefficients from $\mathbb{F}$}. An absolute value $\norm{\cdot}$ on
$\mathcal{L}$ can be defined by setting
\begin{displaymath}
  \norm{\sum_{i=-n}^\infty a_{-i} X^{-i}} = k^n, \quad \norm{0} = 0. 
\end{displaymath}
Under the induced metric, $d(x,y) = \norm{x-y}$, the space
$(\mathcal{L}, d)$ is a complete metric space. Furthermore the
absolute value satisfies for any $x,y \in \mathcal {L}$,
\begin{subequations}
  \begin{equation}
    \label{eq:2}
    \norm{x} \geq 0 \text{ and } \norm{x} = 0 \text{ if and only if }
    x = 0, 
  \end{equation}
  \begin{equation}
    \label{eq:3}
    \norm{xy} = \norm{x}\norm{y},
  \end{equation}
  \begin{equation}
    \label{eq:4}
    \norm{x+y} \leq \max(\norm{x}, \norm{y}).
  \end{equation}
\end{subequations}
Property \eqref{eq:4} is known as the non-Archimedean property. In
fact, equality holds in \eqref{eq:4} whenever $\norm{x} \neq
\norm{y}$. 

As we will be working in finite dimensional vector spaces over
$\mathcal{L}$, we need an appropriate extension of the one-dimensional
absolute value. 
\begin{definition}
  Let $h \in \mathbb{N}$. For any $\mathbf{x} = (x_1, \dots, x_h) \in 
  \mathcal{L}^h$, we define the \emph{height of $x$} to be
  \begin{displaymath}
    \infnorm{\mathbf{x}} = \max\{\norm{x_1}, \dots, \norm{x_h}\}.
  \end{displaymath}
\end{definition}
It is straightforward to see that \eqref{eq:2} and \eqref{eq:4} hold
for \infnorm{\cdot{}}. Of course, when $h=1$, this is the usual
absolute value, and as in the one-dimensional case, \infnorm{\cdot{}}
induces a metric on $\mathcal{L}^h$. When we speak of balls in any of
the spaces $\mathcal{L}^h$, we will mean balls in this metric.

An important consequence of \eqref{eq:4} is that if $C_1$ and $C_2$
are balls in some space $\mathcal{L}^h$, then either $C_1 \cap C_2 =
\emptyset$, $C_1 \subseteq C_2$ or $C_2 \subseteq C_1$. We will refer
to this property as \emph{the ball intersection property.}

In $\mathcal{L}$, the polynomial ring $\mathbb{F}[X]$ plays a r\^{o}le 
analogous to the one played by the integers in the field of real
numbers. Thus, we define \emph{the polynomial part} of a non-zero
element by 
\begin{displaymath}
  \left[\sum_{i=-n}^\infty a_{-i} X^{-i}\right] = \sum_{i=-n}^0 a_{-i}
  X^{-i} 
\end{displaymath}
whenever $n \geq 0$. When $n < 0$, the polynomial part is equal to
zero. Likewise, the polynomial part of the zero element is itself
equal to zero. We define the set
\begin{displaymath}
  I = \left\{ x \in \mathcal{L} : [x]=0 \right\} = \left\{ x \in
  \mathcal{L} : \norm{x} < 1 \right\},
\end{displaymath}
the unit ball in $\mathcal{L}$.

With the above definitions, it makes sense to define the distance to
the polynomial lattice from a point $\mathbf{x} \in \mathcal{L}^h$:
\begin{equation}
  \label{eq:5}
  \intdist{\mathbf{x}} = \min_{\mathbf{p} \in \mathbb{F}[X]^h}
  \infnorm{\mathbf{x}-\mathbf{p}}.
\end{equation}
Since we will be concerned with matrices, we let $m,n \in \mathbb{N}$
be fixed throughout the paper. In the rest of the paper we will need a
number of unspecified constants which may depend on $m$ and $n$. To
avoid cumbersome notation, for such constants, we will only specify
the dependence on parameters other than $m$ and $n$. 

We identify the $m \times n$-matrices with coefficients from
$\mathcal{L}$ with $\mathcal{L}^{mn}$ in the usual way. Matrix
products and inner products are defined as in the real case. Matrices
will be denoted by capital letters, whereas vectors will be denoted by
bold face letters.

In this paper, we are concerned with the Hausdorff dimension (defined
below) of the set of badly approximable systems of linear forms over
$\mathcal{L}$, defined as follows.
\begin{definition}
  The set of matrices
  \begin{displaymath}
    \mathfrak{B}(m,n) = \left\{A \in \mathcal{L}^{mn} : \exists K > 0
      \: \forall \mathbf{q} \in \mathbb{F}[X]^m \setminus \{0\} \: 
      \intdist{\mathbf{q}A}^n > \dfrac{K}{\infnorm{\mathbf{q}}^m}
    \right\}
  \end{displaymath}
  is called \emph{the set of badly approximable elements in
  $\mathcal{L}^{mn}$.}
\end{definition}
On taking $n$'th roots on either side of the defining inequality, we
see that the exponent of $\infnorm{\mathbf{q}}$ on the right hand side
becomes $m/n$. This is exactly the critical exponent in the Laurent
series analogue of the Khintchine--Groshev theorem \cite[Theorem
1]{kristensen}. It is natural to suspect that an analogue of
Dirichlet's theorem exists. This is left as an exercise for the
interested reader.

Let $\mu$ denote the Haar measure on $\mathcal{L}^{mn}$. It is an easy
consequence of \cite[Theorem 1]{kristensen} that $\mathfrak{B}(m,n)$
is a null-set, \emph{i.e.}, $\mu(\mathfrak{B}(m,n))=0$, for any $m,n
\in \mathbb{N}$. This raises the natural question of the Hausdorff
dimension of $\mathfrak{B}(m,n)$, which is shown to be maximal
(Theorem \ref{thm:jarnik} below), thus proving an analogue of
Schmidt's Theorem on badly approximable systems of linear forms over
the real numbers \cite{MR40:1344}.  Niederreiter and Vielhaber
\cite{MR99f:94012} proved using continued fractions that
$\mathfrak{B}(1,1)$ has Hausdorff dimension $1$, \emph{i.e.}, a formal
power series analogue of Jarn\'ik's Theorem
\cite{jarnik28:_zuer_theor_approx}.  The $p$-adic analogue of
Jarn\'ik's Theorem was proven by Abercrombie \cite{MR96g:11078}.

Hausdorff dimension in this setting is defined as follows: Let $E
\subseteq \mathcal{L}^{mn}$. For any countable cover $\mathcal{C}$ of
$E$ with balls $B_i=B(\mathbf{c}_i, \rho_i)$, we define the $s$-length
of $\mathcal{C}$ as the sum
\begin{displaymath}
  l^s(\mathcal{C}) = \sum_{B_i \in \mathcal{C}} \rho_i^s
\end{displaymath}
for any $s \geq 0$. Restricting to covers $\mathcal{C}_\delta$, such that
for some $\delta > 0$, $\rho_i < \delta$ for all $B_i \in
\mathcal{C}_\delta$, we can define an outer measure  
\begin{displaymath}
  \mathcal{H}^s(E) = \lim_{\delta \rightarrow 0} \inf_{\text{covers }
  \mathcal{C}_\delta} l^s(\mathcal{C}_\delta),
\end{displaymath}
commonly called the \emph{Hausdorff $s$-measure} of $E$. It is
straightforward to prove that this is indeed an outer measure. Also,
given a set $E \subseteq \mathcal{L}^{mn}$, the Hausdorff $s$-measure of
$E$ is either zero or infinity for all values of $s \geq 0$, except
possibly one. Furthermore, the Hausdorff $s$-measure of a set is a
non-increasing function of $s$. We define the Hausdorff dimension
$\hdim(E)$ of a set $E \subseteq \mathcal{L}^{mn}$ by
\begin{displaymath}
  \hdim(E) = \inf \left\{ s \geq 0 : \mathcal{H}^s(E) = 0 \right\}.
\end{displaymath}
As in the real case, it can be shown that $\hdim(E) \leq mn$ for any
$E \subseteq \mathcal{L}^{mn}$.

With these definitions, we prove
\begin{theorem}
  \label{thm:jarnik}
  Let $m,n \in \mathbb{N}$. Then,
  \begin{displaymath}
    \text{\emph{dim}}_{\text{\emph{H}}} (\mathfrak{B}(m,n)) = mn.
  \end{displaymath}
\end{theorem}
We will use the method developed by Schmidt \cite{MR33:3793} to
prove the analogous one-dimensional real result, namely the so-called
$(\alpha, \beta)$-games. Schmidt \cite{MR40:1344} subsequently used
this method to prove the multi-dimensional real analogue of Theorem 
\ref{thm:jarnik}. 

The rest of the paper is organised as follows. In section
\ref{sec:notat-defin-prel}, we define $(\alpha, \beta)$-games and some
related concepts and state some results due to Mahler \cite{MR2:350c}
from the appropriate analogue of the geometry of numbers in the
present setting.

The $(\alpha, \beta)$-game has two players, White and Black, with
parameters $\alpha$ and $\beta$ respectively. When played, the game
terminates after infinitely many moves, in a single point in the space
$\mathcal{L}^{mn}$. We prove in section \ref{sec:winn-dimens-badlym},
that for $\alpha$ small enough, player White may ensure that the point
in which the game terminates is an element of $\mathfrak{B}(m,n)$. The
fundamental tools in this proof are a transference principle and a
reduction of the statement to a game which terminates after a finite
number of moves. The transference principle allows us to use the
approximation properties of a matrix to study the approximation
properties of the transpose of the same matrix. The finite game allows
us to show that player White may ensure that all the undesirable
points with $\infnorm{\mathbf{q}}$ less than an appropriate bound can
be avoided. This is the most extensive part of the paper, and the
proof is quite technical.

Finally, in section \ref{sec:hausd-dimens-badlym}, we use the property
from section \ref{sec:winn-dimens-badlym} to show that the dimension
of $\mathfrak{B}(m,n)$ must be greater than or equal to $mn$. Together
with the above remarks, this implies Theorem \ref{thm:jarnik}.

\section{Notation, definitions and preliminary results}
\label{sec:notat-defin-prel}

We now define $(\alpha, \beta)$-games, which will be our main tool in
the proof of Theorem \ref{thm:jarnik}. Let $\Omega = \mathcal{L}^{mn}
\times \mathbb{R}_{\geq 0}$. We call $\Omega$ \emph{the space of
  formal balls in $\mathcal{L}^{mn}$}, where $\omega =
(\mathbf{c},\rho) \in \Omega$ is said to have \emph{centre
  $\mathbf{c}$} and \emph{radius $\rho$}. We define the map $\psi$
from $\Omega$ to the subsets of $\mathcal{L}^{mn}$, assigning a real
closed $\infnorm{\cdot{}}$-ball to the abstract one defined above.
That is, for $\omega = (\mathbf{c},\rho) \in \Omega$,
\begin{displaymath}
  \psi(\omega) = \left\{\mathbf{x} \in \mathcal{L}^{mn} \vert
    \infnorm{\mathbf{x}-\mathbf{c}} \leq \rho\right\}.
\end{displaymath}

\begin{definition}
  Let $B_1, B_2 \in \Omega$. We say that $B_1 =
  (\mathbf{c}_1, \rho_1) \subseteq B_2 = (\mathbf{c}_2, \rho_2)$ if
  $\rho_1 + \infnorm{\mathbf{c}_1 - \mathbf{c}_2} \leq \rho_2$.
\end{definition}
Note that if $B_1 \subseteq B_2$ in $\Omega$, then $\psi(B_1)
\subseteq \psi(B_2)$ as subsets of $\mathcal{L}^{mn}$. Also, we define
for every $\gamma \in (0,1)$ and $B \in \Omega$:
\begin{displaymath}
  B^\gamma = \left\{B' \subseteq B \vert \rho(B') = \gamma
  \rho(B)\right\},
\end{displaymath}
where $\rho(B)$ is the radius of $B$.  We now define the following
game.
\begin{definition}
  \label{defn:alphabetagame}
  Let $S \subseteq \mathcal{L}^{mn}$, and let $\alpha, \beta \in (0,1)$. Let
  Black and White be two players. The \emph{$(\alpha, \beta;
  S)$-game} is played as follows:
  \begin{itemize}
  \item Black chooses a ball $B_1 \in \Omega$.
  \item White chooses a ball $W_1 \in B_1^\alpha$.
  \item Black chooses a ball $B_2 \in W_1^\beta$.
  \item And so on ad infinitum.
  \end{itemize}
  Finally, let $B_i^* = \psi(B_i)$ and $W_i^* =  \psi (W_i)$. If
  $\bigcap_{i=1}^\infty B_i^* = \bigcap_{i=1}^\infty W_i^* \subseteq
  S$, then White wins the game. Otherwise, Black wins the game.
\end{definition}

Our game can be understood in the following way. Initially, Black
chooses a closed ball with radius $\rho_1$. Then, White chooses a ball
with radius $\alpha \rho_1$ inside the first one. Now, Black chooses a
ball with radius $\beta \alpha \rho_1$ inside the one chosen by White,
and so on. In the end, the intersection of these balls will be
non-empty by a simple corollary of Baire's Category Theorem. White
wins the game if this intersection is a subset of $S$. Otherwise,
Black wins.

Because of the unusual topology of $\mathcal{L}^{mn}$, we may
construct distinct elements $(\mathbf{c},\rho), (\mathbf{c}',\rho')
\in \Omega$ such that the corresponding balls in $\mathcal{L}^{mn}$
are the same, \emph{i.e.}, so that $\psi((\mathbf{c},\rho)) =
\psi((\mathbf{c}',\rho'))$ so that the map $\psi$ is not injective.
However, we will often need to consider both the set
$\psi((\mathbf{c},\rho))$ and the formal ball $(\mathbf{c},\rho)$ and
will by abuse of notation denote both by
\begin{displaymath}
  \left\{\mathbf{x} \in \mathcal{L}^{mn} :
    \infnorm{\mathbf{x}-\mathbf{c}} \leq \rho \right\},  
\end{displaymath}
where $\mathbf{c}$ and $\rho$ are understood to be fixed, although
changing these quantities could well have no effect on the set.

The sets of particular interest to us, are sets $S$ such that White
can always win the $(\alpha, \beta; S)$-game. 
\begin{definition}
  A set $S \subseteq \mathcal{L}^{mn}$ is said to be \emph{$(\alpha,
    \beta)$-winning} if White can always win the $(\alpha,
  \beta;S)$-game. $S$ is said to be \emph{$\alpha$-winning} if $S$ is
  $(\alpha,\beta)$-winning for any $\beta \in (0,1)$.
\end{definition}
It is a fairly straightforward matter to see that if $S$ is
$\alpha$-winning for some $\alpha$ and $\alpha' \in (0,\alpha]$, then
$S$ is $\alpha'$-winning. Hence, we may define the maximal $\alpha$
for which a set is $\alpha$-winning.
\begin{definition}
  Let $S \subseteq \mathcal{L}^{mn}$ and let $S^* =\left\{\alpha \in (0,1)
    : S \text{ is $\alpha$-winning}\right\}$. The \emph{winning
    dimension of $S$} is defined as
  \begin{displaymath}
    \windim S =
    \begin{cases}
      0 & \text{if } S^* = \emptyset, \\
      \sup S^* & \text{otherwise.}
    \end{cases}
  \end{displaymath}
\end{definition}

We will first prove that the winning dimension of $\mathfrak{B}(m,n)$
is strictly positive. This will subsequently be used to deduce that
the Hausdorff dimension of $\mathfrak{B}(m,n)$ is maximal. In order to
do this, we study inequalities defined by slightly different matrices.
For any $A \in \mathcal{L}^{mn}$, we define the matrices
\begin{displaymath}
  \matA{A} =
  \begin{pmatrix}
    A & I_m \\ I_n & 0
  \end{pmatrix},
  \qquad \matA{A^*} = 
  \begin{pmatrix}
    A^T & I_n \\ I_m & 0
  \end{pmatrix},
\end{displaymath}
where $I_m$ and $I_n$ denotes the $m \times m$ and $n \times n$
identity matrices respectively. Let \column{j}{A} denote the $j$'th
column of the matrix $A$. In what follows, $\mathbf{q}$ will denote a
vector in $\mathbb{F}[X]^{m+n}$ with coordinates $\mathbf{q} = (q_1,
\dots, q_{m+n})$. Note that $A \in \mathfrak{B}(m,n)$ if and only if there
exists a $K > 0$ such that
\begin{equation}
  \label{eq:8}
  \max_{1 \leq j \leq n} \left(\norm{\mathbf{q} \cdot
  \column{j}{\matA{A}}}\right)^n > \dfrac{K}{\max_{1 \leq i \leq m}
  \left(\norm{q_i}\right)^m}
\end{equation}
for any point in the polynomial lattice $\mathbf{q} \in
\mathbb{F}[X]^{m+n}$ such that the first $m$ coordinates of $\mathbf{q}$
are not all equal to zero.

These matrix inequalities allow us to examine the set $\mathfrak{B}(m,n)$ in
terms of parallelepipeds in $\mathcal{L}^{m+n}$, \emph{i.e.}, sets
defined by inequalities
\begin{equation}
  \label{eq:9}
  \norm{\left(\mathbf{x}A\right)_i} < c_i, \qquad A \in
  \mathcal{L}^{(m+n)^2}, \ c_i > 0, \ i = 1, \dots, m+n,
\end{equation}
where $A$ is invertible and $(\mathbf{x}A)_i$ denotes the $i$'th
coordinate of the vector $\mathbf{x}A$. Inspired by the theory of the
geometry of numbers, we define distance functions
\begin{equation}
  \label{eq:10}
  F_A (\mathbf{x}) := \max_{1 \leq j \leq m+n} \dfrac{1}{c_j}
  \norm{\sum_{i=1}^{m+n} x_i a_{ij}}.
\end{equation}
Also, for any $\lambda > 0$, we define the sets
\begin{displaymath}
  P_A(\lambda) = \left\{ \mathbf{x} \in \mathcal{L}^{m+n} :
    F_A(\mathbf{x}) < \lambda \right\}.
\end{displaymath}
Clearly, $P_A(1)$ is the set defined by \eqref{eq:9}. Also, for
$\lambda' < \lambda$, $P_A(\lambda') \subseteq
P_A(\lambda)$. 

In the setting of the real numbers, distance functions $F_A$ and sets
$P_A$ are studied in the geometry of numbers (see \cite{MR97i:11074}
for an excellent account). For vector spaces over the field of Laurent
series this theory was extensively developed by Mahler in
\cite{MR2:350c}. We will only need a few elementary results, which we
summarise here.

\begin{definition}
  \label{defn:succmin}
  Let $A \in \mathcal{L}^{(m+n)^2}$ be invertible. We define \emph{the
  $j$'th successive minimum $\lambda_j$ of $F_A$} to be
  \begin{multline*}
    \lambda_j = \inf\big\{ \lambda > 0 : P_A(\lambda) \text{ contains $j$
      linearly} \\
    \text{independent } \mathbf{a}_1, \dots, \mathbf{a}_j \in
    \mathbb{F}[X]^{m+n} \big\}.
  \end{multline*}
\end{definition}

We have the following lemma which is a corollary to the result in
\cite[Page 489]{MR2:350c}: 
\begin{lemma}
  \label{lem:sucmin1}
  For any invertible $A \in \mathcal{L}^{(m+n)^2}$,
  \begin{equation}
    \label{eq:12}
    0 < \lambda_1 \leq \dots \leq \lambda_{m+n}.
  \end{equation}
  Furthermore,
  \begin{equation}
    \label{eq:13}
    \lambda_1 \cdots \lambda_{m+n} = \mu(P_A(1))^{-1}.
  \end{equation}
\end{lemma}
It should be noted that Mahler constructs the Haar measure in a
different way from Sprind\v zuk's construction \cite{MR39:6833} used
in \cite{kristensen}. However, as the Haar measure is unique up to a
scaling factor, and since the measure of the unit
$\infnorm{\cdot}$-ball is equal to $1$ in both constructions, the
measures obtained in the two constructions must coincide.

We will need one additional result from \cite[Page 489]{MR2:350c},
relating the successive minima of a parallelepiped to those of its
so-called polar body. 
\begin{lemma}
  \label{lem:sucmin2}
  Let $A \in \mathcal{L}^{(m+n)^2}$ be invertible, let $\lambda_1,
  \dots, \lambda_{m+n}$ denote the successive minima of $F_A$ and let
  $\sigma_1, \dots, \sigma_{m+n}$ denote the successive minima of the
  distance function $F_{A}^*$ defined by
  \begin{displaymath}
    F_A^* (\mathbf{y}) = \sup_{\mathbf{x} \neq 0}
    \dfrac{\norm{\mathbf{x} \cdot \mathbf{y}}}{F_A(\mathbf{x})}.
  \end{displaymath}
  Then,
  \begin{displaymath}
    \lambda_m \sigma_{n+1} = 1.
  \end{displaymath}
\end{lemma}
The definition of a polar body can be taken to be the one implicit in
the statement of Lemma \ref{lem:sucmin2}.

\section{The winning dimension of $\mathfrak{B}(m,n)$}
\label{sec:winn-dimens-badlym}

In this section, we will prove that the winning dimension of
$\mathfrak{B}(m,n)$ is strictly positive. We will obtain an explicit
lower bound on the winning dimension.  For the rest of this section,
let $n,m \in \mathbb{N}$ be fixed and $\alpha, \beta \in (0,1)$ be
such that $\gamma = k^{-1} + \alpha\beta - (k^{-1} + 1)\alpha > 0$.

We now begin the game. Black starts by choosing a ball $B_1$ of radius
$\rho = \rho(B_1)$. Clearly the set $B_1$ is bounded, so we may fix a
$\sigma > 0$ such that for all $A \in B_1$, $\infnorm{A} \leq \sigma$.
We will construct a strategy for player White depending on a constant
$R > R_0(\alpha,\beta,\rho,\sigma) \geq 1$, which we will choose
later. We use subsequently
\begin{displaymath}
  \delta = R^{-m(m+n)^2}, \quad \delta^* = R^{-n(m+n)^2}, \quad
  \tau = \dfrac{m}{m+n}.
\end{displaymath}
Let $B_k, B_h \subseteq \mathcal{L}^{mn}$  be balls occurring in the
$(\alpha, \beta)$-game chosen by Black such that $\rho(B_k) <
R^{-(m+n)(\tau +i)}$ and $\rho(B_h) < R^{-(m+n)(1 +j)}$ for
some $i,j \in \mathbb{N}$. We will show that White can play in such a
way that the following properties hold for $i,j \in \mathbb{N}$:
\begin{itemize}
\item \label{item:1} For $A \in B_k$, there are no $\mathbf{q} \in
  \mathbb{F}[X]^{m+n}$ such that the inequalities 
  \begin{subequations}
    \begin{equation}
      \label{eq:14}
      0 < \max_{1 \leq l \leq m} \left\{\norm{q_l}\right\} < \delta
      R^{n(\tau+i)}
    \end{equation}
    and
    \begin{equation}
      \label{eq:15}
      \max_{1 \leq l' \leq n} \left\{\norm{\mathbf{q} \cdot
      \column{l'}{\matA{A}}}\right\} < \delta R^{-m(\tau +i) - n}
    \end{equation}
  \end{subequations}
  both hold.
\item \label{item:2} For $A \in B_h$, there are no $\mathbf{q} \in
  \mathbb{F}[X]^{m+n}$ such that the inequalities
  \begin{subequations}
    \begin{equation}
      \label{eq:16}
      0 < \max_{1 \leq l' \leq n} \left\{\norm{q_{l'}}\right\} <
      \delta^* R^{m(1+j)}
    \end{equation}
    and
    \begin{equation}
      \label{eq:17}
      \max_{1 \leq l \leq m} \left\{\norm{\mathbf{q} \cdot
      \column{l}{\matA{A^*}}}\right\} < \delta^* R^{-n(1 +j) -
      m}.
    \end{equation}
  \end{subequations}
  both hold.
\end{itemize}
If White follows a strategy such that \eqref{eq:14} and \eqref{eq:15}
are avoided for all $i \in \mathbb{N}$, she will win the $(\alpha,
\beta; \mathfrak{B}(m,n))$-game. Indeed, given a $\mathbf{q} \in
\mathbb{F}[X]^{m+n}$ with the first $m$ coordinates, $q_1, \dots, q_m$
say, not all equal to zero, we can find an $i \in \mathbb{N}$ such
that
\begin{equation}
  \label{eq:11}
  \delta R^{n(\tau+i-1)} \leq \max_{1 \leq l \leq m}
  \left\{\norm{q_l}\right\} < \delta R^{n(\tau+i)}.
\end{equation}
This immediately implies that \eqref{eq:14} holds for this $i$, so
that \eqref{eq:15} must be false. Hence, by \eqref{eq:11},
\begin{alignat*}{2}
  \max_{1 \leq l' \leq n} \left\{\norm{\mathbf{q} \cdot
  \column{l'}{\matA{A}}}\right\}^n &\geq \dfrac{ \delta^{m+n}
  R^{-mn(\tau +i) - n^2 + mn(\tau +i) -mn}}{\max_{1 \leq l \leq m}
  \left\{\norm{q_l}\right\}^m} \\
  &\geq \dfrac{\delta^{m+n} R^{-n^2 -mn}}{\max_{1 \leq l
  \leq m} \left\{\norm{q_l}\right\}^m} \\
  & > \dfrac{K}{\max_{1 \leq l \leq m} \left\{\norm{q_l}\right\}^m}
\end{alignat*}
for any $K \in (0, \delta^{m+n}R^{-n^2-mn})$, the matrix $A$ is in
$\mathfrak{B}(m,n)$ by \eqref{eq:8}.

For the remainder of this section, we will construct a strategy for
White ensuring that \eqref{eq:14} and \eqref{eq:15}
(resp.~\eqref{eq:16} and \eqref{eq:17}) cannot hold for any $i$
(resp.~$j$). We define for any $i \in \mathbb{N}$:
\begin{itemize}
\item $B_{k_i}$ to be the first ball chosen by Black with
  $\rho(B_{k_i}) < R^{-(m+n)(\tau +i)}$.
\item $B_{h_i}$ to be the first ball chosen by Black with
  $\rho(B_{h_i}) < R^{-(m+n)(1 +i)}$.
\end{itemize}
Since $\tau < 1$, these balls occur such that $B_{k_0}
\supseteq B_{h_0} \supseteq B_{k_1} \supseteq B_{h_1} \supseteq
\cdots$. By choosing $R$ large enough, we can ensure that the
inclusions are proper.

Since
\begin{displaymath}
  \delta R^{n \tau} = R^{-m(m+n)^2 + nm(m+n)^{-1}} =
  R^{-m\left((m+n)^2 - \frac{n}{m+n}\right)} < 1,
\end{displaymath}
\eqref{eq:14} has no solutions for $i = 0$. Hence, White can certainly
play in such a way \eqref{eq:14} and \eqref{eq:15} have no polynomial
solutions when $A \in B_{k_0}$. We will construct White's strategy in
such a way that:
\begin{enumerate}
\item \label{item:3} Given the beginning of a game $B_1 \supseteq W_1
  \supseteq \cdots \supseteq B_{k_0} \supseteq \cdots \supseteq
  B_{k_i}$ such that \eqref{eq:14} and \eqref{eq:15} have no polynomial
  solutions for any $A \in B_{k_i}$, White can play in such a way that
  \eqref{eq:16} and \eqref{eq:17} have no polynomial solutions for any $A
  \in B_{h_i}$.
\item \label{item:4} Given the beginning of a game $B_1 \supseteq W_1
  \supseteq \cdots \supseteq B_{k_0} \supseteq \cdots \supseteq
  B_{h_i}$ such that \eqref{eq:16} and \eqref{eq:17} have no polynomial
  solutions for any $A \in B_{h_i}$, White can play in such a way that
  \eqref{eq:14} and \eqref{eq:15} have no polynomial solutions for any $A
  \in B_{k_{i+1}}$.
\end{enumerate}

Our first lemma guarantees that we need only consider solutions to the
equations in certain subspaces of $\mathcal{L}^{m+n}$.
\begin{lemma}
  \label{lem:subspace1}
  Let $B_1 \supseteq W_1 \supseteq \cdots \supseteq B_{k_i}$ be the
  start of a game such that \eqref{eq:14} and \eqref{eq:15} have no
  polynomial solutions for any $A \in B_{k_i}$. The set
  \begin{displaymath}
    \left\{\mathbf{q} \in \mathbb{F}[X]^{m+n} : \eqref{eq:16} \text{ and } 
    \eqref{eq:17} \text{ hold for } j = i \text{ for some } A \in
    B_{k_i} \right\}
  \end{displaymath}
  contains at most $m$ linearly independent points.
\end{lemma}

\begin{proof}
  Assume that there are linearly independent $\mathbf{q}_1, \dots,
  \mathbf{q}_{m+1} \in \mathbb{F}[X]^{m+n}$ such that \eqref{eq:16}
  and \eqref{eq:17} hold for $A_1, \dots, A_{m+1} \in B_{k_i}$.  The
  absolute value of the first $n$ coordinates must be less than
  $\delta^* R^{m(1 + i)}$ by \eqref{eq:16}, and as $\infnorm{A_u} \leq
  \sigma$ for $u=1, \dots, m+1$, \eqref{eq:17} and the structure of
  $\matA{A_u^*}$ guarantee that there is a constant $K_1(\sigma) > 0$
  such that
  \begin{equation}
    \label{eq:18}
    \infnorm{\mathbf{q}_u} \leq K_1(\sigma) \delta^* R^{m(1 +
      i)} \quad \text{for } 1 \leq u \leq m+1.
  \end{equation}
  Let $C$ be the centre of $B_{k_i}$. For any $A \in B_{k_i}$,
  \begin{equation}
    \label{eq:19}
    \infnorm{\column{l}{\matA{A^*}} - \column{l}{\matA{C^*}}} \leq
    \rho(B_{k_i}) < R^{-(m+n)(\tau + i)} \quad \text{for } 1 \leq l
    \leq m.
  \end{equation}
  Now, as \eqref{eq:17} holds for the vectors, \eqref{eq:18} and
  \eqref{eq:19} imply that for $u = 1, \dots, m+1$,
  \begin{multline}
    \label{eq:20}
    \max_{1 \leq l \leq m} \bigg\{\bigg\Vert \mathbf{q}_u \cdot
    \column{l}{\matA{C^*}}  \bigg\Vert\bigg\} \\
    \leq \max_{1 \leq l \leq m} \left\{\norm{\mathbf{q}_u \cdot
        \column{l}{\matA{A_u^*}}}, 
      \norm{\mathbf{q}_u \cdot \left( \column{l}{\matA{C^*}} -
          \column{l}{\matA{A_u^*}}\right)}\right\} \\
    \leq \max \left\{ \delta^* R^{-n(1 +i) -m}, K_1(\sigma)
      \delta^* R^{m(1 +i)} R^{- (m+n)(\tau +i)} \right\} \\
    \leq K_2(\sigma) \delta^* R^{-n(1 + i)},
  \end{multline}
  where $K_2(\sigma) > 0$. If needed, we may increase the right hand
  side, so that without loss of generality, $K_2(\sigma) > 1$.

  We define the parallelepiped
  \begin{multline*}
    P = \bigg\{\mathbf{y} \in \mathcal{L}^{m+n} : \max_{1 \leq l' \leq n}
      \left\{\norm{y_{l'}}\right\} < R^{m(1 +i)}, \\
      \max_{1 \leq l \leq m} \left\{\norm{\mathbf{y} \cdot
          \column{l}{\matA{C^*}}}\right\} < R^{-n(1 +i)}
    \bigg\}, 
  \end{multline*}
  along with the corresponding distance function $F_C$ and the
  successive minima $\lambda_1, \dots, \lambda_{m+n}$. By
  \eqref{eq:20}, $\lambda_{m+1} \leq K_2(\sigma) \delta^*$. For $n=1$,
  $0 < \lambda_{m+1} \leq K_2(\sigma) R^{-(m+1)^2}$, which by Lemma
  \ref{lem:sucmin1} gives a contradiction by choosing $R$ large
  enough.

  Hence, we may assume that $n > 1$. Let
  \begin{multline*}
    P^* = \bigg\{\mathbf{x} \in \mathcal{L}^{m+n} : \max_{1 \leq l \leq
      m} \left\{\norm{x_{l}}\right\} < R^{n(1 +i)},\\
    \max_{1 \leq l' \leq n} \left\{\norm{\mathbf{x} \cdot
        \column{l'}{\matA{C}}}\right\} < R^{-m(1 +i)}
    \bigg\}. 
  \end{multline*}
  This set admits the distance function $F_C^*$ defined in Lemma
  \ref{lem:sucmin2} as the two bodies, $P$ and $P^*$, are mutually
  polar (see \cite{MR2:350c}).

  Let $\sigma_1, \dots, \sigma_{m+n}$ denote the successive minima of
  $P^*$. By Lemma \ref{lem:sucmin1} and Lemma \ref{lem:sucmin2},
  \begin{multline*}
    \sigma_1 \leq \left(\sigma_1 \cdots
      \sigma_{n-1}\right)^{\frac{1}{n-1}} = \mu(P^*)^{\frac{-1}{n-1}}
    \left(\sigma_n \cdots \sigma_{m+n}\right)^{\frac{-1}{n-1}}\\
    \leq \mu(P^*)^{\frac{-1}{n-1}} \sigma_n^{-\frac{m+1}{n-1}} =
    \mu(P^*)^{\frac{-1}{n-1}} \lambda_{m+1}^{\frac{m+1}{n-1}} \leq
    \mu(P^*)^{\frac{-1}{n-1}}
    \left(K_2(\sigma) \delta^*\right)^{\frac{m+1}{n-1}} \\
    \leq K_3(\sigma) R^{-(m+n)^2(m+1)} = K_3(\sigma) \delta
    R^{-(m+n)^2},
  \end{multline*}
  where $K_3(\sigma) > 0$. Hence, there is a $\mathbf{q} \in
  \mathbb{F}[X]^{m+n}\setminus \{0\}$ with 
  \begin{displaymath}
    \max_{1 \leq l \leq m} \left\{\norm{q_{l}}\right\} < K_3(\sigma)
    \delta R^{-(m+n)^2} R^{n(1 +i)}
  \end{displaymath}
  and
  \begin{equation}
    \label{eq:31}
    \max_{1 \leq l' \leq n} \left\{\norm{\mathbf{q} \cdot
    \column{l'}{\matA{C}}}\right\} < K_3(\sigma) \delta R^{-(m+n)^2} 
    R^{-m(1 +i)} < 1,
  \end{equation}
  when we choose $R$ large enough. But \eqref{eq:31} implies that
  $\max_{1 \leq l \leq m} \left\{\norm{q_{l}}\right\} > 0$, since
  otherwise the last $n$ coordinates would also be equal to $0$,
  whence $\mathbf{q} = 0$. This gives a contradiction, as we have
  found a solution to \eqref{eq:14} and \eqref{eq:15}.
\end{proof}

In a completely analogous way, we can prove:
\begin{lemma}
  \label{lem:subspace2}
  Let $B_1 \supseteq W_1 \supseteq \cdots \supseteq B_{h_i}$ be the
  start of a game such that \eqref{eq:16} and \eqref{eq:17} have no
  polynomial solutions for any $A \in B_{h_i}$. The set
  \begin{multline*}
    \big\{\mathbf{q} \in \mathbb{F}[X]^{m+n} : \eqref{eq:14} \text{ and }
    \eqref{eq:15}\\ \text{ hold with $i$ replaced by $i+1$ for some } A
    \in B_{h_i} \big\}
  \end{multline*}
  contains at most $n$ linearly independent points.
\end{lemma}

We will now reduce the statement that White has a strategy such that
Step \ref{item:3} on page \pageref{item:3} is possible, to the
statement that White can win a certain finite game. The converse Step
\ref{item:4} is analogous. 

Once again, we assume that $B_1 \supseteq W_1 \supseteq \cdots
\supseteq B_{k_i}$ is the beginning of a game such that we have
avoided polynomial solutions to all relevant inequalities so far. Now,
it is sufficient for White to avoid solutions $\mathbf{q} \in
\mathbb{F}[X]^{m+n}$ to \eqref{eq:16} and \eqref{eq:17} with
\begin{displaymath}
  \delta^* R^{m(1 +i -1)} \leq \max_{1 \leq l' \leq n}
  \left\{\norm{q_{l'}} \right\} < \delta^* R^{m(1 +i)},
\end{displaymath}
as solutions have been avoided for all vectors $\mathbf{q}$ with
\begin{displaymath}
  \max_{1 \leq l' \leq n} \left\{\norm{q_{l'}} \right\} < \delta^* 
  R^{m(1 +i -1)}
\end{displaymath}
in the preceeding steps by assumption.  Hence we need
only consider $\mathbf{q} \in \mathbb{F}[X]^{m+n}$ for which
\begin{equation}
  \label{eq:38}
  \delta^* R^{m(1 +i -1)} \leq \infnorm{\mathbf{q}}.
\end{equation}

By Lemma \ref{lem:subspace1}, the set of $\mathbf{q}$ satisfying
\eqref{eq:16} and \eqref{eq:17} is contained in some $m$-dimensional
subspace. Let $\{\mathbf{y}_1, \dots, \mathbf{y}_m\}$ be an
orthonormal basis for this space and write all $\mathbf{q}$ in this
subspace satisfying \eqref{eq:38} in the form $\mathbf{q} = t_1
\mathbf{y}_1 + \cdots + t_m \mathbf{y}_m$, $t_1,\dots, t_m \in
\mathcal{L}$. Immediately,
\begin{equation}
  \label{eq:21}
  \delta^* R^{m(1 +i -1)} \leq \max_{1 \leq l' \leq m}
  \left\{\norm{t_{l'}}\right\}.
\end{equation}
White needs to avoid solutions to the inequalities
\begin{equation}
  \label{eq:23}
  \max_{1 \leq l \leq m}\norm{\sum_{l'=1}^m t_{l'}
    \left(\mathbf{y}_{l'} \cdot \column{l}{\matA{A^*}}\right)} <
  \delta^* R^{-n(1 +i) -m}. 
\end{equation}
This matrix inequality may be solved using Cramer's Rule \cite[Chapter
XIII, Theorem 4.4]{MR86j:00003}. This theorem shows that \eqref{eq:23}
is soluble if for $l'=1, \dots, m$,
\begin{displaymath}
  \norm{t_{l'}}\norm{D} = \norm{t_{l'} D} \leq \delta^* R^{-n(1
    +i) -m} \max_{1 \leq l \leq m} \left\{\norm{D_{l, l'}}\right\},
\end{displaymath}
where $D$ denotes the determinant of the matrix with entries
$\mathbf{y}_{l'} \cdot \column{l}{\matA{A^*}}$ and $D_{l,l'}$ denotes
the $(l,l')$'th co-factor of this determinant. By \eqref{eq:21}, it is
sufficient to avoid
\begin{multline}
  \label{eq:22}
  \norm{D} \leq R^{-n(1 +i) - m - m(1 +i - 1)} \max_{1 
    \leq l,l' \leq m} \left\{\norm{D_{l,l'}}\right\} \\
  = R^{-(m+n)(1 +i)}  \max_{1 \leq l,l' \leq m}
  \left\{\norm{D_{l,l'}}\right\}. 
\end{multline}
We define the following finite game:
\begin{definition}
  \label{defn:finitegame}
  Let $\mathbf{y}_1, \dots, \mathbf{y}_m \in \mathcal{L}^{m+n}$ be a
  set of orthonormal vectors. Let $B \subseteq \mathcal{L}^{mn}$ be a
  ball with $\rho(B) < 1$ such that for any $A \in B$, $\infnorm{A}
  \leq \sigma$. Let $\mu > 0$ and let $\alpha, \beta \in (0,1)$ with
  $k^{-1} + \alpha\beta -(k^{-1}+1)\alpha > 0$. White and Black take
  turns according to the rules of the game in Definition
  \ref{defn:alphabetagame}, choosing balls inside $B$, but the game
  terminates when $\rho(B_t) < \mu \rho(B)$. White wins the game if
  \begin{displaymath}
    \norm{D} > \rho(B) \mu \max_{1 \leq l,l' \leq m}
    \left\{\norm{D_{l,l'}}\right\} 
  \end{displaymath}
  for any $A \in B_t$.
\end{definition}

If White can win the game in Definition \ref{defn:finitegame} for any
$\mu \in (0, \mu^*)$ for some $\mu^* = \mu^*(\alpha, \beta, \sigma) >
0$, then White can guarantee that \eqref{eq:22} does not hold for any
$A \in B_{h_i}$. To see this, let $B = B_{k_i}$ and let
\begin{displaymath}
  \mu = \dfrac{R^{-(m+n)(1 +i)}}{\rho(B)} \leq (\alpha
  \beta)^{-1} R^{-n}.
\end{displaymath}
Choosing $R$ large enough, this will be less than $\mu^*$. It remains
to be shown, that such a $\mu^*$ exists. We will do this by induction.

Let $A\in \mathcal{L}^{mn}$, $v \in \{1, \dots,m\}$ and $\{\mathbf{y}_1, 
\dots,\mathbf{y}_m\}$ be the orthonormal system from Definition
\ref{defn:finitegame}. By considering all possible choices of 
$1\leq i_1 < \cdots < i_v \leq m$ and $1 \leq j_1 < \cdots < j_v \leq
m$, we obtain 
$\binom{m}{v}^2$ matrices 
\begin{equation}
  \label{eq:25}
  \begin{pmatrix}
    \mathbf{y}_{i_1} \cdot \column{j_1}{\matA{A^*}} & \cdots &
    \mathbf{y}_{i_1} \cdot \column{j_v}{\matA{A^*}} \\
    \vdots & & \vdots \\
    \mathbf{y}_{i_v} \cdot \column{j_1}{\matA{A^*}} & \cdots &
    \mathbf{y}_{i_v} \cdot \column{j_v}{\matA{A^*}}
  \end{pmatrix}
\end{equation}
For each $v \in \{1, \dots, m\}$, we define the function $M_v : 
\mathcal{L}^{mn} \rightarrow \mathcal{L}^{\binom{m}{v}^2}$ to have as it's
coordinates the determinants of the matrices in \eqref{eq:25} in some
arbitrary but fixed order. Furthermore, define
\begin{displaymath}
  M_{-1}(A) = M_0(A) = (1),
\end{displaymath}
the standard unit vector in $\mathcal{L}^{\binom{m}{0}^2} =
\mathcal{L}$.  For $K \subseteq \mathcal{L}^{mn}$, we define
\begin{displaymath}
  M_v(K) = \max_{A \in K} \infnorm{M_v(A)}.
\end{displaymath}

We will prove a series of lemmas, culminating in a proof that under
appropriate conditions, player White may always win the game in
Definition \ref{defn:finitegame} (Lemma \ref{lem:winfinite1} below). 

In the following, assume that $v > 0$ and that there exists a
$\mu_{v-1}$ such that
\begin{equation}
  \label{eq:26}
  \infnorm{M_{v-1}(A)} > \rho (B) \mu_{v-1} M_{v-2}(B_{i_{v-1}}) 
\end{equation}
for all $A \in B_{i_{v-1}}$ for an appropriate $B_{i_{v-1}}$
occurring in the game. 

\begin{lemma}
  \label{lem:winfinite2}
  Let $\epsilon >0$ and let $B' \subseteq B_{i_{v-1}}$ be a ball
  of radius 
  \begin{displaymath}
    \rho(B') < \epsilon \mu_{v-1} \rho(B_{i_{v-1}}). 
  \end{displaymath}
  Then 
  \begin{displaymath}
    \infnorm{M_{v-1}(A) - M_{v-1}(A')} < \epsilon \rho(B) \mu_{v-1}
    M_{v-2}(B_{i_{v-1}}) 
  \end{displaymath}
  for any $A, A' \in B'$.
\end{lemma}

\begin{proof}
  Consider first for a fixed $A \in B'$ and a fixed $x \in
  \mathcal{L}$ the quantity
  \begin{displaymath}
    \infnorm{M_{v-1}(A + x E_{ij}) - M_{v-1}(A)},
  \end{displaymath}
  where $E_{ij}$ denotes the matrix with $1$ in the $ij$'th
  entry and zeros elsewhere. On considering an individual
  coordinate of the vector $M_{v-1}(A + x E_{ij}) - M_{v-1}(A)$ and
  applying the ultra-metric inequality \eqref{eq:4}, it is seen that
  \begin{equation}
    \label{eq:6}
    \infnorm{M_{v-1}(A + x E_{ij}) - M_{v-1}(A)} \leq
    \norm{x} M_{v-2}(B_{i_{v-1}}). 
  \end{equation}
  The factor $M_{v-2}(B_{i_{v-1}})$ is an upper bound on the co-factor
  corresponding to the $ij$'th minor. When $\norm{x} = 1$, These
  quantities are discrete analogues of the partial derivatives of
  $M_{v-1}$, and the upper bound \eqref{eq:6} implies that the
  function does not vary wildly.
  
  We may pass from one $A$ matrix of $B'$ to another $A'$ by changing
  one coordinate at a time, \emph{i.e.}, by performing a string of
  $mn$ operations $A \mapsto A+(A'_{ij} - A_{ij}) E_{ij}$.  Using
  these operations, we define a finite sequence of matrices by
  $A^{(1,1)} = A+(A'_{11} - A_{11}) E_{11}$, $A^{(2,1)} =
  A^{(1,1)}+(A'_{21} - A^{(1,1)}_{21}) E_{21}$ and so on, so that
  $A^{(m,n)} = A'$. We now obtain,
  \begin{multline*}
    \infnorm{M_{v-1}(A) - M_{v-1}(A')} = \bigg\Vert M_{v-1}(A) -
    M_{v-1}(A^{(1,1)}) + M_{v-1}(A^{(1,1)}) \\ - \dots -
    M_{v-1}(A^{((m-1),n)}) + M_{v-1}(A^{((m-1),n)}) -
    M_{v-1}(A')\bigg\Vert_\infty 
  \end{multline*}
  Here, each matrix in the arguments of $M_{v-1}$ differ from the
  preceding one in only one place.  Applying \eqref{eq:6} and the
  ultra-metric inequality \eqref{eq:4} $mn$ times,
  \begin{multline*}
    \infnorm{M_{v-1}(A) - M_{v-1}(A')} \leq \infnorm{A - A'}
    M_{v-2}(B_{i_{v-1}}) \\ 
    < \epsilon \mu_{v-1} \rho(B_{i_{v-1}}) M_{v-2}(B_{i_{v-1}}) \leq
    \epsilon \rho(B) \mu_{v-1} M_{v-2}(B_{i_{v-1}}).
  \end{multline*}
\end{proof}

\begin{corollary}
  \label{cor:finite1}
  For a ball $B' \subseteq B_{i_{v-1}}$ with radius $\rho(B') < 
  \tfrac{1}{2} \mu_{v-1} \rho(B_{i_{v-1}})$, we have
  \begin{displaymath}
    \infnorm{M_{v-1}(A')} > \tfrac{1}{2} M_{v-1}(B')
  \end{displaymath}
  for any $A' \in B'$.
\end{corollary}

\begin{proof}
  Apply Lemma \ref{lem:winfinite2} with $\epsilon = \tfrac{1}{2}$
  and use \eqref{eq:26}.
\end{proof}
  
Now, we define 
\begin{displaymath}
  D_v(A) = \det 
  \begin{pmatrix}
    \mathbf{y}_1 \cdot \column{1}{\matA{A^*}} & \cdots &
    \mathbf{y}_1 \cdot \column{v}{\matA{A^*}} \\
    \vdots && \vdots \\
    \mathbf{y}_v \cdot \column{1}{\matA{A^*}} & \cdots &
    \mathbf{y}_v \cdot \column{v}{\matA{A^*}}
  \end{pmatrix}.
\end{displaymath}
Clearly, this is a function of the $nv$ variables $a_{11}, \dots,
a_{n1}, \dots, a_{nv}$. We define the \emph{discrete gradient of
  $D_v$} to be the vector
\begin{displaymath}
  \nabla D_v(A) = 
  \begin{pmatrix}
    D_v(A+E_{11}) - D_v(A) \\
    \vdots \\
    D_v(A+E_{mn}) - D_v(A),
  \end{pmatrix} \in \mathcal{L}^{mn},
\end{displaymath}
where $E_{ij} \in \mathcal{L}^{mn}$ denotes the matrix having $1$ as
the $ij$'th entry and zeros elsewhere.
\begin{corollary}
  \label{cor:finite2}
  With $B'$ as in Lemma \ref{lem:winfinite2} and $A', A'' \in B'$,
  we have
  \begin{displaymath}
    \infnorm{\nabla D_v(A') - \nabla D_v (A'')} \leq K_4
    \infnorm{M_{v-1} (A') - M_{v-1}(A'')}
  \end{displaymath}
  for some $K_4 > 0$ depending only on $m$ and $n$.
\end{corollary}

\begin{proof}
  Note, that the coordinates of $\nabla D_v(A)$ are linear
  combinations of the coordinates of $M_{v-1}(A)$ for any $A$.
\end{proof}

The discrete gradient turns out to be the key ingredient in the
proof. We will need the following lemma: 
\begin{lemma}
  \label{lem:gradpositiv}
  Let $B' \subseteq B_{i_{v-1}}$ be a ball such that 
  \begin{equation}
    \label{eq:27}
    \rho(B') < \tfrac{1}{2} \mu_{v-1} \rho(B_{i_{v-1}}).
  \end{equation}
  Let $A' \in B'$ be such that
  \begin{equation}
    \label{eq:28}
    \infnorm{M_v(A')} < \tfrac{1}{8} M_{v-1} (B').
  \end{equation}
  Furthermore, assume that the maximum $\infnorm{M_{v-1}(A')}$ is
  attained by the absolute value of the coordinate which is the
  determinant 
  \begin{displaymath}
    d_v = \det
    \begin{pmatrix}
      \mathbf{y}_1 \cdot \column{1}{\matA{{A'}^*}} & \cdots &
      \mathbf{y}_1 \cdot \column{v-1}{\matA{{A'}^*}} \\
      \vdots && \vdots \\
      \mathbf{y}_{v-1} \cdot \column{1}{\matA{{A'}^*}} & \cdots &
      \mathbf{y}_{v-1} \cdot \column{v-1}{\matA{{A'}^*}} 
    \end{pmatrix}
    .
  \end{displaymath}
  Then
  \begin{displaymath}
    \infnorm{\nabla D_v(A')} > K_5(\sigma) M_{v-1} (B')
  \end{displaymath}
  for some $K_5(\sigma) > 0$.
\end{lemma}
  
\begin{proof}
  Let $\mathbf{z} \in \mathcal{L}^{m+n}$. Consider the following
  quantity 
  \begin{multline}
    \label{eq:29}
    \Phi(\mathbf{z}) = \left\Vert \det
      \begin{pmatrix}
        \mathbf{y}_1 \cdot \column{1}{\matA{{A'}^*}} & \cdots &
        \mathbf{y}_1 \cdot \left(\column{v}{\matA{{A'}^*}} +
          \mathbf{z}\right)\\ 
        \vdots && \vdots \\
        \mathbf{y}_{v} \cdot \column{1}{\matA{{A'}^*}} & \cdots &
        \mathbf{y}_{v} \cdot \left(\column{v}{\matA{{A'}^*}}
          +\mathbf{z}\right)
      \end{pmatrix}
    \right. \\
    - \left. \det 
      \begin{pmatrix}
        \mathbf{y}_1 \cdot \column{1}{\matA{{A'}^*}} & \cdots &
        \mathbf{y}_1 \cdot \column{v}{\matA{{A'}^*}} \\
        \vdots && \vdots \\
        \mathbf{y}_{v} \cdot \column{1}{\matA{{A'}^*}} & \cdots &
        \mathbf{y}_{v} \cdot \column{v}{\matA{{A'}^*}} 
      \end{pmatrix}
    \right\Vert \\
    = \norm{\left(\sum_{h=1}^v (-1)^{h+1} d_h \mathbf{y}_h \right)
      \cdot \mathbf{z}}. 
  \end{multline}    
  The last equality follows on expanding the determinants in the
  last column, where the $d_i$ are taken from the coordinates of
  $M_{v-1}(A')$ and $d_v$ is the special coordinate for which the
  maximum absolute value is attained.
  
  Let $\mathbf{z}_1 = d_1 \mathbf{y}_1 + \cdots + d_{v-1}
  \mathbf{y}_{v-1} + X d_v \mathbf{y}_v$, where $X \in \mathcal{L}$ is
  the power series consisting solely of the indeterminate $X$. We have
  assumed that $\norm{d_i} \leq \norm{d_v} < k \norm{d_v} = \norm{X
    d_v}$ for all $i = 1, \dots, v-1$.  Hence, since the
  $\mathbf{y}_i$ are assumed to be orthonormal,
  \begin{multline*}
    \norm{\left(\sum_{h=1}^v (-1)^{h+1} d_h \mathbf{y}_h \right)
      \cdot \mathbf{z}_1}\\
    = \norm{\left(\sum_{h=1}^v (-1)^{h+1} d_h \mathbf{y}_h \right)
      \cdot \left(d_1 \mathbf{y}_1 + \cdots + X d_v
        \mathbf{y}_v\right)} \\ 
    = k \norm{d_v}^2 = k \infnorm{M_{v-1}(A')}^2 > \tfrac{k}{4}
    M_{v-1}(B')^2 
  \end{multline*}
  by Corollary \ref{cor:finite1}.
  
  We wish to interpret $\Phi(\mathbf{z})$ as a discrete analogue of
  the directional derivative along a vector in $\mathcal{L}^{mn}$.
  Furthermore, we need to obtain a lower bound on this quantity for
  some direction. In order to be able to make this interpretation, we
  need to find a lower bound on $\Phi(\mathbf{z}_2)$, where
  $\mathbf{z}_2$ is of the form $(z_1, \dots, z_n, 0, \dots, 0)$,
  \emph{i.e.}, where the last $m$ coordinates are zero. For such
  vectors, considering the difference in \eqref{eq:29} corresponds to
  considering the difference $\norm{D_v(A + \hat{Z}_2) - D_v(A)}$,
  where $\hat{Z}_2 \in \mathcal{L}^{mn}$ is the matrix which has the
  vector $(z_1, \dots, z_n)$ as its $v$'th row and zeros elsewhere, so
  that the matrix $A + \hat{Z}_2 \in \mathcal{L}^{mn}$ is the matrix
  $A$ with the entries of the $v$'th row shifted by the first $n$
  coordinates of $\mathbf{Z}_2$.  When $\infnorm{\hat{Z}_2} = 1$, this
  quantity is exactly the discrete partial derivative of $D_v$ in
  direction $\hat{Z}_2$ evaluated at $A$.
  
  Because of the special form of the $\column{l}{\matA{{A'}^*}}$, we 
  may write
  \begin{displaymath}
    \mathbf{y}_h = \mathbf{y}_h^0 + \lambda_{h1}
    \column{1}{\matA{{A'}^*}} + \cdots + \lambda_{hm}
    \column{m}{\matA{{A'}^*}}, 
  \end{displaymath}
  where the $\mathbf{y}_h^0$ have zeros on the last $m$ coordinates.
  Since the $\mathbf{y}_h$ are orthonormal, certainly for all $h,l$,
  $\norm{\lambda_{h,l}} \leq 1$. Also, there is a constant
  $K_5(\sigma) > 0$ such that $\infnorm{\mathbf{y}_h^0} \leq
  \tfrac{1}{8} K_5(\sigma)^{-1}$ for all $h$. We define
  $\mathbf{z}_2 = d_1 \mathbf{y}_1^0 + \cdots + d_{v-1}
  \mathbf{y}_{v-1}^0 + X d_v \mathbf{y}_v^0$, which clearly has the
  required form.
  
  Now, 
  \begin{multline}
    \label{eq:30}
    \Phi(\mathbf{z}_2) = \norm{\left(\sum_{h=1}^v (-1)^{h+1} d_h
        \mathbf{y}_h \right) \cdot \mathbf{z}_2}\\
    = \norm{\left(\sum_{h=1}^v (-1)^{h+1} d_h \mathbf{y}_h \right)
      \cdot (\mathbf{z}_2 - \mathbf{z}_1 + \mathbf{z}_1)} \\
    \geq \tfrac{k}{4} M_{v-1}(B')^2 - \norm{\left(\sum_{h=1}^v
        (-1)^{h+1} d_h \mathbf{y}_h \right) \cdot (\mathbf{z}_1 - 
      \mathbf{z}_2)}. 
  \end{multline}
  In order to produce a good lower bound for $\Phi(\mathbf{z}_2)$,
  we will produce a good upper bound on the last term of the above.
  We know that
  \begin{displaymath}
    \mathbf{z}_1 - \mathbf{z}_2 = \sum_{l=1}^m \left(\sum_{h=1}^v
      d'_h \lambda_{hl}\right) \column{l}{\matA{{A'}^*}},
  \end{displaymath}
  where $d'_h = d_h$ for $h = 1,\dots, v-1$ and $d_v =
  Xd_v$. Furthermore for $l=l, \dots, m$, by simple calculation,
  \begin{multline*}
    \norm{\left(\sum_{h=1}^v (-1)^{h+1} d_h \mathbf{y}_h \right) 
      \cdot \column{l}{\matA{{A'}^*}}} \\
    = \norm{\det
      \begin{pmatrix}
        \mathbf{y}_1 \cdot \column{1}{\matA{{A'}^*}} & \cdots &
        \mathbf{y}_1 \cdot \column{v-1}{\matA{{A'}^*}} &
        \mathbf{y}_1 \cdot \column{l}{\matA{{A'}^*}} \\
        \vdots && \vdots \\
        \mathbf{y}_{v} \cdot \column{1}{\matA{{A'}^*}} & \cdots &
        \mathbf{y}_{v} \cdot \column{v-1}{\matA{{A'}^*}} &
        \mathbf{y}_{v} \cdot \column{l}{\matA{{A'}^*}} 
      \end{pmatrix}
      } \\
    \leq \infnorm{M_v(A')} < \tfrac{1}{8} M_{v-1}(B') 
  \end{multline*}
  by choice of $A'$. Hence, as $\norm{d'_h} \leq k M_{v-1}(B')$,
  \begin{multline*}
    \norm{\left(\sum_{h=1}^v (-1)^{h+1} d_h \mathbf{y}_h \right)
      \cdot (\mathbf{z}_1 - \mathbf{z}_2)} \\
    = \norm{\sum_{l=1}^m \left(\sum_{h'=1}^v d'_{h'}
        \lambda_{h'l}\right)  
      \left(\sum_{h=1}^v (-1)^{h+1} d_h \mathbf{y}_h \right) \cdot 
      \column{l}{\matA{{A'}^*}}} < \tfrac{k}{8} M_{v-1}(B')^2.
  \end{multline*}
  Together with \eqref{eq:30} this implies
  \begin{equation}
    \label{eq:43}
    \Phi(\mathbf{z}_2) > \tfrac{k}{8} M_{v-1}(B')^2.
  \end{equation}
  
  We wish to use the discrete directional derivative to obtain a lower
  bound on the discrete gradient. Let $\mathbf{z} \in
  \mathcal{L}^{m+n}$ be some vector of the form $(z_1, \dots, z_n, 0,
  \dots, 0)$, so that $\mathbf{z}$ corresponds to a matrix $\hat{Z}_v
  \in \mathcal{L}^{mn}$ with $(z_1, \dots, z_n)$ as its $v$'th row and
  zeros elsewhere.  Suppose further that $\infnorm{\mathbf{z}} =
  \infnorm{\hat{Z}_v} = 1$. It is simple to show that $\infnorm{\nabla
    D_v(A')} \geq \Phi(\mathbf{z})$.
  
  Let $\log_k$ denote the logarithm to base $k$. We normalise
  $\mathbf{z}_2$ by $X^{-\log_k \infnorm{\mathbf{z}_2}}$, where $X$ is
  again the indeterminate in the power series expansions. In this way,
  we obtain a vector in $\mathcal{L}^{m+n}$ corresponding to a matrix
  $\hat{Z}_2 \in \mathcal{L}^{mn}$ with $\infnorm{\hat{Z}_2} = 1$.
  Now, note that by \eqref{eq:29}, for any $x \in \mathcal{L}$ and any
  $\mathbf{z} \in \mathcal{L}^{m+n}$, $\Phi(x \mathbf{z}) = \norm{x}
  \Phi(\mathbf{z})$. But as $\infnorm{\mathbf{z}_2} \leq \tfrac{1}{8}
  K_5(\sigma)^{-1} k M_{v-1}(B')$, we get by \eqref{eq:43}
  \begin{multline*}
    \infnorm{\nabla D_v(A)} \geq \Phi(\mathbf{z}_2 \, X^{-\log_k 
      \infnorm{\mathbf{z}_2}}) =
    \dfrac{\Phi(\mathbf{z}_2)}{\infnorm{\mathbf{z}_2}}\\
    > \dfrac{\tfrac{k}{8} 
      M_{v-1}(B')^2}{\tfrac{1}{8} K_5(\sigma)^{-1} k M_{v-1}(B')} = 
    K_5(\sigma) M_{v-1}(B').
  \end{multline*}
  This completes the proof.
\end{proof}

We are now ready to prove that player White can win the game defined
in Definition \ref{defn:finitegame}.

\begin{lemma}
  \label{lem:winfinite1}
  Let $\{\mathbf{y}_1, \dots, \mathbf{y}_m\} \subseteq
  \mathcal{L}^{m+n}$ be a set of orthonormal vectors. Let $B \subseteq
  \mathcal{L}^{mn}$ be a ball, $\rho(B) = \rho_0 < 1$, such that for
  some $\sigma > 0$, $\infnorm{A} < \sigma$ for any $A \in B$. Let
  $\alpha, \beta \in (0,1)$ with $k^{-1} + \alpha\beta -(k^{-1} + 1)
  \alpha > 0$. Assume that $0 \leq v \leq m$.
  
  There exists a $\mu_v = \mu_v(\alpha,\beta,\sigma) \in (0,2]$
  for which White can play the game in Definition
  \ref{defn:finitegame} in such a way that for the first ball
  $B_{i_v}$ with $\rho(B_{i_v}) < \rho_0 \mu_v$,
  \begin{displaymath}
    \infnorm{M_v(A)} > \rho_0 \mu_v M_{v-1}(B_{i_v})
  \end{displaymath}
  for any $A \in B_{i_v}$.
\end{lemma}

The slightly cumbersome notation $B_{i_v}$ is used in order to make
the connection with \eqref{eq:26}, which we will use in the proof,
explicit. Of course, the additional subscript plays no r\^{o}le in the
statement of the lemma.

\begin{proof}
  We will prove the lemma by induction. Clearly, the lemma holds for
  $v=0$. Hence, we use \eqref{eq:26} as our induction hypothesis, and
  so we have the above results as our disposal.

  Recall that $\gamma = k^{-1} + \alpha\beta - (k^{-1}+1)\alpha > 0$
  and let
  \begin{displaymath}
    \epsilon = \dfrac{\gamma}{8}\dfrac{K_5(\sigma)}{K_4} >
    0. 
  \end{displaymath}
  Furthermore, let
  \begin{displaymath}
    i_v = \min \left\{ i \in \mathbb{N} : i > i_{v-1}, \rho(B_i) <
    \min(\tfrac{1}{2}, \epsilon) \mu_{v-1}
    \rho(B_{i_{v-1}})\right\}. 
  \end{displaymath}
  By appropriately choosing a constant $K_6(\alpha,\beta,\sigma) >
  0$, we have
  \begin{equation}
    \label{eq:36}
    \rho(B_{i_v}) \geq K_6(\alpha,\beta,\sigma)\rho_0.
  \end{equation}

  Using the induction hypothesis, Corollary \ref{cor:finite2} and
  Lemma \ref{lem:winfinite2}, for any $A', A'' \in B_{i_v}$, we have 
  \begin{equation}
    \label{eq:41}
    \infnorm{\nabla D_v (A') - \nabla D_v(A'')} < K_4 \epsilon
    \rho_0 \mu_{v-1} M_{v-2}(B_{i_{v-1}})
    < \tfrac{\gamma}{8} K_5(\sigma) M_{v-1}(B_{i_v}).
  \end{equation}
  We now let 
  \begin{displaymath}
    \mu_v = \min\left\{\tfrac{1}{8}, \tfrac{\gamma}{8} \alpha\beta
      K_6(\alpha,\beta, \sigma), \tfrac{3\gamma}{8} K_5(\sigma)
      K_6(\alpha,\beta, \sigma)K_7(\sigma) \right\} > 0, 
  \end{displaymath} 
  where $K_7(\sigma)>0$ is to be chosen later. Assume that there
  exists an $A' \in B_{i_v}$ for which the assertion of the lemma does
  not hold. That is,
  \begin{displaymath}
    \infnorm{M_v(A')} \leq \rho_0 \mu_v M_{v-1}(B_{i_v}).
  \end{displaymath}
  In this case, we will prove that White has a strategy which will
  eliminate such elements in a finite number of moves.
  
  By choice of $i_v$, \eqref{eq:27} holds. Since $\rho_0 < 1$,
  \eqref{eq:28} holds. By rearranging the $\mathbf{y}_i$, we can
  without loss of generality assume that the condition on the
  determinant in Lemma \ref{lem:gradpositiv} holds. Hence,
  \begin{equation}
    \label{eq:32}
    \infnorm{\nabla D_v (A')} > K_5(\sigma) M_{v-1}(B_{i_v}).
  \end{equation}
  
  Let $\nabla' = \nabla D_v(A')$, and let $D_i$ and $C_i$ denote the
  centres of $W_i$ and $B_i$ respectively. White can play in such a
  way that 
  \begin{equation}
    \label{eq:33}
    \norm{(C_i - D_i) \cdot \nabla'} \geq k^{-1}(1-\alpha) \rho(B_i)
    \infnorm{\nabla'}.
  \end{equation}
  Indeed, there are points $D_i \in B_i$ with $\infnorm{C_i - D_i}
  \geq k^{-1}(1- \alpha) \rho(B_i)$ and such that $B(D_i, \alpha
  \rho(B_i)) \subseteq B_i$. This guarantees \eqref{eq:33}. Also, no
  matter how Black plays
  \begin{equation}
    \label{eq:34}
    \norm{(C_{i+1} - D_i) \cdot \nabla'} \leq (1-\beta) \rho(W_i)
    \infnorm{\nabla'},
  \end{equation}
  since Black cannot choose the next centre further away from
  $D_i$. Hence,
  \begin{multline*}
    \norm{(C_{i+1} - C_i) \cdot \nabla'} \geq \left(k^{-1}(1-\alpha) -
      \alpha(1-\beta)\right) \rho(B_i) \infnorm{\nabla'}\\
    = \gamma \rho(B_i) \infnorm{\nabla'} > 0.
  \end{multline*}
  We choose $t_0 \in \mathbb{N}$ such that $\alpha\beta
  \tfrac{\gamma}{2} < (\alpha \beta)^{t_0} \leq \tfrac{\gamma}{2}$.
  Player White can ensure that
  \begin{equation}
    \label{eq:35}
    \norm{(C_{i+t_0} - C_i) \cdot \nabla'} \geq \gamma \rho(B_i)
    \infnorm{\nabla'} > 0.
  \end{equation}
  This follows from \eqref{eq:33}, \eqref{eq:34} and the fact that
  $\gamma > 0$ so that player White can ensure that the bound in
  \eqref{eq:34} is preserved for the next $t_0$ steps. White will play
  according to such a strategy. Furthermore $\rho(B_{i+t_0}) \leq
  \tfrac{\gamma}{2}\rho(B_i)$, so for any $A \in B_{i+t_0}$,
  \begin{multline}
    \label{eq:37}
    \norm{(A - C_i) \cdot \nabla'} \geq \norm{(C_{i+t_0} - C_i) \cdot
      \nabla'} - \norm{(A - C_{i+t_0}) \cdot \nabla'}\\
    \geq
    \tfrac{\gamma}{2} \rho(B_i) \infnorm{\nabla'}.
  \end{multline}

  Now, for any $A \in B_{i_v}$, 
  \begin{multline}
    \label{eq:7}
    \norm{\left(A - C_{i_v}\right) \cdot \nabla D_v(A)} \leq
    \infnorm{A - C_{i_v}} \norm{\nabla D_v(A)} \\
    \leq \rho(B_{i_v}) \max_{\substack{1 \leq l \leq m \\ 1 \leq l'
        \leq n}} \left\{\norm{D_v\left(A + E_{ll'}\right)},
      \norm{D_v(A)}\right\} \leq K_7(\sigma) \norm{D_v(A)}
  \end{multline}  
  for some $K_7(\sigma) > 0$. Also, 
  \begin{multline}
    \label{eq:40}
    \norm{(A - C_{i_v}) \cdot \nabla'} = \norm{(A - C_{i_v}) \cdot
      \nabla D_v(A) + (A - C_{i_v}) \cdot (\nabla' - \nabla D_v(A))}
    \\
    \leq \max\{\norm{(A - C_{i_v}) \cdot \nabla D_v(A)}, \norm{(A -
      C_{i_v}) \cdot (\nabla' - \nabla D_v(A))} \} \\
    \leq \norm{(A - C_{i_v}) \cdot \nabla D_v(A)} + \norm{(A -
      C_{i_v}) \cdot (\nabla D_v(A) - \nabla')}
  \end{multline}
  
  Combining inequalities \eqref{eq:7} and \eqref{eq:40}, we obtain for
  some $K_7(\sigma) > 0$,
  \begin{equation}
    \label{eq:42}
    \norm{D_v(A)} \geq K_7(\sigma)\big(\norm{(A - C_{i_v})\cdot
    \nabla'} - \norm{(A-C_{i_v}) \cdot (\nabla D_v(A) -
    \nabla')}\big).
  \end{equation}
  Now by \eqref{eq:37} and \eqref{eq:32},
  \begin{equation}
    \label{eq:39}
    \norm{(A - C_{i_v})\cdot \nabla'} \geq \tfrac{\gamma}{2} \rho(B_{i_v})
    \infnorm{\nabla'} \geq \tfrac{\gamma}{2} \rho(B_{i_v}) K_5(\sigma)
    M_{v-1}(B_{i_v}).
  \end{equation}
  By \eqref{eq:41},
  \begin{multline*}
    \norm{(A-C_{i_v}) \cdot (\nabla D_v(A) -  \nabla')} \leq
    \rho(B_{i_v}) \infnorm{\nabla D_v(A) - \nabla D_v(A')} \\
    \leq \rho(B_{i_v}) \tfrac{\gamma}{8} K_5 (\sigma)
    M_{v-1}(B_{i_v}). 
  \end{multline*}
  Combining this with \eqref{eq:42} and \eqref{eq:39}, 
  \begin{displaymath}
    \norm{D_v(A)} \geq K_7(\sigma) \tfrac{3 \gamma}{8}
    \rho(B_{i_v}) K_5 (\sigma) M_{v-1}(B_{i_v}).
  \end{displaymath}
  Inserting \eqref{eq:36} into this expression, we find that 
  \begin{displaymath}
    \norm{D_v(A)} \geq \tfrac{3 \gamma}{8} K_5 (\sigma)
    K_6(\alpha,\beta,\sigma) K_7(\sigma) \rho_0 M_{v-1}(B_{i_v})
    \geq \mu_v \rho_0 M_{v-1}(B_{i_v}),
  \end{displaymath}
  by choice of $\mu_v$. This completes the proof of Lemma
  \ref{lem:winfinite1}.
\end{proof}

Note that Lemma \ref{lem:winfinite1} immediately implies:
\begin{theorem}
  \label{thm:windim}
  Let $\alpha, \beta \in (0,1)$ with $k^{-1} + \alpha\beta - (k^{-1} +
  1)\alpha > 0$ and let $m,n \in \mathbb{N}$. White can win the game
  in Definition \ref{defn:finitegame} and hence the $(\alpha, \beta;
  \mathfrak{B}(m,n))$-game. In particular,
  \begin{displaymath}
    \windim\left( \mathfrak{B}(m,n) \right) \geq \dfrac{1}{k+1}.
  \end{displaymath}
\end{theorem}
\begin{proof}
  The first part follows from Lemma \ref{lem:winfinite1} with $v=m$
  and the obvious analogue for the other step in the strategy for the
  $(\alpha, \beta; \mathfrak{B}(m,n))$-game. The lower bound on the
  winning dimension follows as $k^{-1} + \alpha\beta - (k^{-1} +
  1)\alpha > 0$ for any $\beta \in (0,1)$ and any $\alpha < 1/(k+1)$.
\end{proof}

\section{The Hausdorff dimension of $\mathfrak{B}(m,n)$}
\label{sec:hausd-dimens-badlym}
  
In this final section, we will prove that if $\alpha > 0$, then any
$\alpha$-winning set in $\mathcal{L}^{mn}$ has full Hausdorff
dimension.  By Theorem \ref{thm:windim}, this will imply Theorem
\ref{thm:jarnik}. To do this, we change our viewpoint to that of
player Black. We will for each step of the game examine a number of
different possible directions for the game under the assumption that
player White is following a winning strategy. This will give rise to a
particularly rich subset of the $\alpha$-winning set for which we may
estimate the Hausdorff dimension.

\begin{theorem}
  \label{thm:hdim}
  Let $\beta \in (0,1)$ and let $N(\beta) \in \mathbb{N}$ be such that
  any ball $B \subseteq \mathcal{L}^{mn}$ of radius $\rho$ contains
  $N(\beta)$ pairwise disjoint balls of radius $\beta \rho$. Let $S
  \subseteq \mathcal{L}^{mn}$ be $(\alpha,\beta)$-winning. Then 
  \begin{displaymath}
    \text{\emph{dim}}_{\text{\emph{H}}}(S) \geq \dfrac{\log
      N(\beta)}{\abs{\log \alpha \beta}}.
  \end{displaymath}
\end{theorem}

\begin{proof}
  Let $\Lambda = \{0,\dots,N(\beta)-1\}^\mathbb{N}$ and let $(i_j) \in
  \Lambda$. For each ball $W_j$ chosen by White, we pick $N(\beta)$
  disjoint balls in $W_j^\beta$ which we enumerate by elements from
  the set $\{0, \dots, N(\beta)-1\}$. We restrict the choice of moves
  for player Black to these $N(\beta)$ possibilities. In this way, we
  obtain for each element $\lambda \in \Lambda$ a point $A(\lambda)
  \in \mathcal{L}^{mn}$. As we may assume that White is following a
  winning strategy, for each $\lambda \in \Lambda$, $A(\lambda) \in
  S$. We will label balls chosen by player Black by the sequence
  leading to them, \emph{i.e.}, $B_l = B(i_1, \dots, i_l)$, where
  $i_1, \dots, i_j \in \{0, \dots, N(\beta) - 1\}$. As distinct
  sequences give rise to disjoint balls from some point in the game
  and onwards, distinct points $\lambda, \lambda' \in \Lambda$ give
  rise to different points $A(\lambda), A(\lambda') \in S$.
  
  Let
  \begin{displaymath}
    S^* = \bigcup_{\lambda \in \Lambda} \{A(\lambda)\} \subseteq S.
  \end{displaymath}
  We define a surjective function $f: S^* \rightarrow [0,1]$ by 
  \begin{displaymath}
    A \mapsto x=0.i_1i_2\dots \quad \text{where } A=A(i_1,i_2,\dots)
  \end{displaymath}
  where $0.i_1 i_2\dots$ is the base $N(\beta)$ expansion of $x$. We
  extend this function to all subsets of $\mathcal{L}^{mn}$ in the
  following way. For $T \subseteq S^*$, let $f(T) = \bigcup_{A\in T}
  f(A)$. For $R\subseteq \mathcal{L}^{mn}$, let $f(R) = f(R\cap S^*)$.
  
  Let $\mathcal{C} = (B_l)_{l \in \mathbb{N}}$ be a cover of $S$ with
  balls, where $B_l$ has radius $\rho_l$. Clearly, $\mathcal{C}^* =
  (B_l \cap S^*)_{l \in \mathbb{N}}$ is a cover of $S^*$. Mapping to
  the interval, we find that $f(\mathcal{C}^*) = (f(B_l \cap S^*))_{l
    \in \mathbb{N}} = (f(B_l))_{l\in \mathbb{N}}$ is a cover of
  $[0,1]$. Thus, the union of the sets $f(B_l)$ has outer Lebesgue
  measure $\ell$ greater than $1$, so by sub-additivity
  \begin{equation}
    \label{eq:24}
    \sum_{l=1}^\infty \ell \left(f(B_l)\right) \geq 1.
  \end{equation}

  Now, let
  \begin{displaymath}
    j_l = \left[\dfrac{\log 2\rho_l}{\log \alpha\beta}\right].
  \end{displaymath}
  For $\rho_l$ sufficiently small, we have $j_l > 0$ and $\rho_l <
  (\alpha\beta)^{j_l}$. Hence, by the ball intersection property,
  $B_l$ is contained in at most one ball of the form $B_l(i_1, \dots,
  i_{j_l})$. But such a ball clearly maps into an interval of length
  $N(\beta)^{-j_l}$. Hence, $\ell(f(B_l)) \leq N(\beta)^{-j_l}$. By
  \eqref{eq:24}, we have
  \begin{multline*}
    1 \leq \sum_{l=1}^\infty \ell\left(f(B_l)\right) \leq
    \sum_{l=1}^\infty N(\beta)^{-j_l} \\
    = \sum_{l=1}^\infty
    N(\beta)^{-\left[\tfrac{\log 2\rho_l}{\log \alpha\beta}\right]}
    \leq 2^{\tfrac{\log N(\beta)}{\abs{\log \alpha \beta}}}
    \sum_{l=1}^\infty \rho_l^{\tfrac{\log N(\beta)}{\abs{\log
          \alpha\beta}}}. 
  \end{multline*}
  Now, for any such cover $\mathcal{C}$ with small enough
  balls, the $s$-length $l^s(\mathcal{C}) >0$ for $s = \tfrac{\log
  N(\beta)}{\abs{\log \alpha\beta}}$. Hence,
  \begin{displaymath}
    \hdim (S) \geq \dfrac{\log N(\beta)}{\abs{\log \alpha\beta}}.
  \end{displaymath}
\end{proof}

Theorem \ref{thm:hdim} allows us to prove that $\hdim
(\mathfrak{B}(m,n)) = mn$ and thus complete the proof of Theorem
\ref{thm:jarnik}: 
\begin{proof}[Proof of Theorem \ref{thm:jarnik}]
  By Theorem \ref{thm:hdim}, we need only estimate the number
  $N(\beta)$ to get a lower bound for the Hausdorff dimension. This is
  a simple combinatorial problem. By scaling and translation, we note
  that it suffices to consider the $\infnorm{\cdot}$-ball $B(0,1) =
  I^{mn}$.
  
  We choose the number $i \in \mathbb{Z}$ such that $k^{i-1} \leq
  \beta < k^i$ and consider the balls $B(c,\beta) \subseteq I^{mn}$
  where $c \in X^{i+1}\mathbb{F}[X]^{mn}$. By choice of $i$ and the
  ball intersection property, these are clearly disjoint. Furthermore,
  counting these balls we see that
  \begin{displaymath}
    N(\beta) = \left(k^{-i-1}\right)^{mn} = \dfrac{1}{k^{mn}}
    \dfrac{1}{\left(k^i\right)^{mn}} \asymp \dfrac{1}{\beta^{mn}}.
  \end{displaymath}
  Hence, by Theorem \ref{thm:hdim},
  \begin{displaymath}
    \hdim(\mathfrak{B}(m,n)) \geq \dfrac{mn\abs{\log \beta}}{\abs{\log
    \alpha} + \abs{\log \beta}} \xrightarrow[\beta \rightarrow 0]{} mn.
  \end{displaymath}
  This completes the proof.
\end{proof}

\subsection*{Acknowledgements}
  I thank my Ph.D.~supervisors R.~Nair and M.~Weber
  for their encouragement. In particular, I thank R.~Nair for
  suggesting the topic of this paper. I also thank M.~M.~Dodson for
  his comments on the manuscript.
  
  The bulk of the research presented in this paper was done while I
  was a student at the University of Liverpool. Many changes and
  corrections to the paper have been made since then.


\begin{thebibliography}{99}

\bibitem{MR96g:11078}
  \textsc{A.~G. Abercrombie}, \textit{Badly approximable $p$-adic
  integers}, Proc. Indian Acad. Sci. Math. Sci. {\bf105~(2)}
  (1995), 123--134. 

\bibitem{MR97i:11074}
  \textsc{J.~W.~S. Cassels}, \textit{An introduction to the geometry of
    numbers}. Springer-Verlag, Berlin, 1997.

\bibitem{jarnik28:_zuer_theor_approx}
  \textsc{V.~Jarn\'ik}, \textit{Zur metrischen {T}heorie der
    diophantischen Approximationen}, Prace Mat.--Fiz. (1928--1929),
  91--106. 

\bibitem{kristensen}
  \textsc{S.~Kristensen}, \textit{On the well-approximable matrices
    over a field of formal series}, Math. Proc. Cambridge
  Philos. Soc. {\bf135~(2)} (2003), 255--268.
  
\bibitem{MR86j:00003}
  \textsc{S.~Lang}, \textit{Algebra {\rm(2nd
      Edition)}}. Addison-Wesley Publishing Co., Reading, Mass., 1984.

\bibitem{MR2:350c}
  \textsc{K.~Mahler}, \textit{An analogue to {M}inkowski's geometry of 
    numbers in a field of series}, Ann. of Math. (2) {\bf42} (1941),
  488--522. 
  
\bibitem{MR99f:94012}
  \textsc{H.~Niederreiter and M.~Vielhaber}, \textit{Linear complexity 
  profiles: {H}ausdorff dimensions for almost perfect profiles and
  measures for general profiles}, J. Complexity {\bf13~(3)} (1997),
  353--383. 

\bibitem{MR40:1344}
  \textsc{W.~M. Schmidt}, \textit{Badly approximable systems of linear 
  forms}, J. Number Theory {\bf1} (1969), 139--154.

\bibitem{MR33:3793}
  \textsc{W.~M. Schmidt}, \textit{On badly approximable numbers and
  certain games}, Trans. Amer. Math. Soc. {\bf123} (1966), 178--199.

\bibitem{MR39:6833}
  \textsc{V.~G. Sprind{\v{z}}uk}, \textit{Mahler's problem in metric 
    number theory}. American  Mathematical Society, Providence, R.I.,
  1969. 
\end{thebibliography}
\end{document}